\documentclass[11pt]{amsart}
\usepackage[T1]{fontenc}
\usepackage[utf8]{inputenc}
\usepackage{amsmath,amssymb,amsthm}
\usepackage[a4paper,margin=1.15in]{geometry}
\usepackage{booktabs}
\usepackage{float}
\usepackage{fancyvrb}
\usepackage{newunicodechar}
\usepackage[colorlinks,linkcolor=blue,citecolor=blue,urlcolor=blue]{hyperref}

\newunicodechar{ℕ}{\ensuremath{\mathbb{N}}}
\newunicodechar{ℤ}{\ensuremath{\mathbb{Z}}}
\newunicodechar{→}{\ensuremath{\to}}
\newunicodechar{↔}{\ensuremath{\leftrightarrow}}
\newunicodechar{∀}{\ensuremath{\forall}}
\newunicodechar{∃}{\ensuremath{\exists}}
\newunicodechar{∈}{\ensuremath{\in}}
\newunicodechar{∨}{\ensuremath{\lor}}
\newunicodechar{∧}{\ensuremath{\land}}
\newunicodechar{≤}{\ensuremath{\le}}
\newunicodechar{∣}{\ensuremath{\mid}}
\newunicodechar{·}{\ensuremath{\cdot}}
\newunicodechar{ε}{\ensuremath{\varepsilon}}
\newunicodechar{δ}{\ensuremath{\delta}}
\fvset{fontsize=\footnotesize,xleftmargin=1em}

\newtheorem{theorem}{Theorem}[section]
\newtheorem{lemma}[theorem]{Lemma}
\newtheorem{proposition}[theorem]{Proposition}
\newtheorem{corollary}[theorem]{Corollary}
\newtheorem{conjecture}[theorem]{Conjecture}
\theoremstyle{remark}
\newtheorem{remark}[theorem]{Remark}
\newtheorem{example}[theorem]{Example}

\newcommand{\rev}{\operatorname{rev}}
\newcommand{\cC}{\mathcal{C}}
\newcommand{\cN}{\mathcal{N}}
\newcommand{\FF}{\mathbb{F}}
\newcommand{\ZZ}{\mathbb{Z}}
\newcommand{\one}{\mathbf{1}}
\newcommand{\tp}{\mathsf{T}}
\newcommand{\seg}[2]{\chi[#1,#2]}

\begin{document}

\title{Nilpotent tridiagonal matrices with entries in $\{-1,0,1\}$}
\author{Marcel Blattner}
\address{Applied AI Research Lab, Lucerne University of Applied Sciences and Arts, Suurstoffi~1, 6343 Rotkreuz, Switzerland}
\email{marcel.blattner@hslu.ch}
\urladdr{https://orcid.org/0009-0001-7188-673X}
\subjclass[2020]{Primary 15B36; Secondary 11A55, 15A18, 15B35, 94A55}
\keywords{Bohemian matrices, nilpotent matrices, tridiagonal matrices, continued fractions, folding lemma, sign patterns}

\begin{abstract}
We determine all nilpotent tridiagonal matrices with entries in $\{-1,0,1\}$. The problem reduces to describing the set $\cN_n$ of pairs $(d,b)\in\{-1,0,1\}^n\times\{-1,1\}^{n-1}$ for which the tridiagonal matrices with diagonal $d$ and with products $b_k$ of opposite off-diagonal entries have characteristic polynomial $\lambda^n$. We show that $\cN_n$ is obtained from $\cN_{\lfloor n/2\rfloor}$ by one of two explicit ``folding'' operations, according to the parity of $n$. As a consequence, $\cN_n$ is non-empty precisely when $n=1$ or $3\cdot2^j-1\le n\le2^{j+2}-1$ for some $j\ge0$, in which case it has $2^{\lfloor\log_2n\rfloor}$ elements. This yields a formula for the number of nilpotent $n\times n$ tridiagonal matrices with entries in $\{-1,0,1\}$, which begins $1,9,73,585,4809,39113$. The case $d=0$ is a conjecture of Corless, established in an earlier paper. The main new ingredient is a rigidity phenomenon: a pair as above, of size $n\ge2$, whose characteristic polynomial is congruent to $\lambda^n$ modulo $4$ is necessarily a fold of a smaller pair. It follows that such a pair is nilpotent as soon as its characteristic polynomial is congruent to $\lambda^n$ modulo $2^{1+\lfloor\log_2n\rfloor}$. Finally we show that the eigenvalue zero of an $m\times m$ matrix in Corless's skew-symmetric family has multiplicity at most $2^{v_2(m+1)}-1$, with equality attained. The classification, the rigidity statement, the cardinality of the sets $\cN_n$ and the criterion modulo $2^{1+\lfloor\log_2n\rfloor}$ have been formalised in Lean~4.
\end{abstract}

\maketitle

\section{Introduction}

\subsection{The classification}
Let $n\ge1$. This paper is concerned with the $n\times n$ tridiagonal matrices
\begin{equation}\label{eq:M}
M=\begin{pmatrix}
d_1 & a_1 & & & \\
c_1 & d_2 & a_2 & & \\
 & c_2 & \ddots & \ddots & \\
 & & \ddots & \ddots & a_{n-1}\\
 & & & c_{n-1} & d_n
\end{pmatrix},\qquad d_k,a_k,c_k\in\{-1,0,1\},
\end{equation}
and with the question of which of them are nilpotent. Families of matrices whose entries are drawn from a fixed finite set, the \emph{population}, are known as \emph{Bohemian} matrices, and have been studied intensively in recent years, largely by experiment; see \cite{Chan} and the references therein, where in particular the problem of understanding the number of nilpotent matrices in such a family is raised.

Our starting point is the work of Corless \cite{Corless} on the skew-symmetric tridiagonal matrices with superdiagonal entries in $\{1,i\}$. He showed that such a matrix can be nilpotent only when its dimension is of the form $2^K-1$, gave a recursive construction of $2^{K-1}$ nilpotent matrices in that dimension, and conjectured that there are no others. This conjecture was established in \cite{Blattner}, with a formalisation in Lean \cite{BlattnerCorless}, and independently by Shen \cite{Shen}. Corless's matrices have the same characteristic polynomials as the matrices \eqref{eq:M} with $d=0$ and all $a_kc_k=\pm1$, and it is natural to ask what happens when the diagonal is permitted to be non-zero. It turns out that the answer has the same recursive flavour, but that a second construction now enters, which doubles the dimension $n$ rather than sending it to $2n+1$. As a result, the dimensions in which unreduced nilpotent matrices exist, that is to say nilpotent matrices \eqref{eq:M} none of whose off-diagonal entries vanishes, are governed by the binary expansion of $n$, in a way that is perhaps not obvious at first sight (Corollary~\ref{cor:sizes}).

We begin with some routine reductions. The characteristic polynomial $\det(\lambda-M)$ of \eqref{eq:M} depends on the off-diagonal entries only through the products $b_k:=a_kc_k$, and it factors into the characteristic polynomials of two smaller matrices of the same type whenever some $b_k$ vanishes. It therefore suffices to understand the case in which all the $b_k$ are non-zero, and then the question of whether $M$ is nilpotent depends only on the pair
\[
(d,b)\in\cC_n:=\{-1,0,1\}^n\times\{-1,1\}^{n-1}.
\]
More generally, for a pair $p=(d,b)\in R^n\times R^{n-1}$ over a commutative ring $R$ we write $\chi_p\in R[\lambda]$ for the characteristic polynomial of the tridiagonal matrices with diagonal $d$ and with $a_kc_k=b_k$; we call $n$ the \emph{size} of $p$, and we call $p$ \emph{nilpotent} if $\chi_p=\lambda^n$. Our task is thus to describe the sets
\[
\cN_n:=\{p\in\cC_n:\chi_p=\lambda^n\}.
\]

There are two ways to produce a nilpotent pair from a smaller one. Write $\rev x$ for the reversal of a finite sequence $x$. Given a pair $p=(d,b)$ of size $n$ and signs $\varepsilon,\delta\in\{-1,1\}$, we define the \emph{odd fold} $\Phi_\varepsilon(p)$ and the \emph{even fold} $\Psi_\delta(p)$ of $p$, of sizes $2n+1$ and $2n$ respectively, by
\begin{align}
\Phi_\varepsilon(p)&:=\bigl((d,\,0,\,\rev d),\ (b,\,\varepsilon,\,-\varepsilon,\,\rev b)\bigr),\label{eq:odd}\\
\Psi_\delta(p)&:=\bigl((d_1,\dots,d_{n-1},\,d_n-\delta,\,d_n+\delta,\,d_{n-1},\dots,d_1),\ (b,\,-1,\,\rev b)\bigr).\label{eq:even}
\end{align}
A short computation (Lemma~\ref{lem:fold}) shows that
\begin{equation}\label{eq:fold}
\chi_{\Phi_\varepsilon(p)}=\lambda\,\chi_p^2\qquad\text{and}\qquad\chi_{\Psi_\delta(p)}=\chi_p^2,
\end{equation}
so that the folds of a nilpotent pair are again nilpotent. The odd fold maps $\cC_n$ into $\cC_{2n+1}$. The even fold of a pair $p\in\cC_n$ lies in $\cC_{2n}$ precisely when the last diagonal entry $d_n$ of $p$ vanishes. Corless's construction is the odd fold, restricted to $d=0$.

Our main result asserts that there are no further nilpotent pairs.

\begin{theorem}[Classification]\label{thm:main}
We have $\cN_1=\{o\}$, where $o:=((0),())$, and for every $h\ge1$,
\begin{align*}
\cN_{2h+1}&=\bigl\{\Phi_\varepsilon(p):\ p\in\cN_h,\ \varepsilon\in\{-1,1\}\bigr\},\\
\cN_{2h}&=\bigl\{\Psi_\delta(p):\ p\in\cN_h\text{ with last diagonal entry }0,\ \delta\in\{-1,1\}\bigr\}.
\end{align*}
Moreover, each element of $\cN_{2h+1}$ and of $\cN_{2h}$ arises in this way from exactly one $p$ and one sign.
\end{theorem}

In other words, the nilpotent pairs of size $n$ are obtained from the $1\times1$ zero matrix by applying one fold for each binary digit of $n$ after the leading one, an odd fold for each digit $1$ and an even fold for each digit $0$, subject to the proviso that an even fold may only be applied to a pair whose last diagonal entry vanishes. For instance, for $n=13=1101_2$ the sizes of the successive pairs are $1,3,6,13$, and no obstruction arises. For $n=10=1010_2$ they are $1,2,5,10$; but the pairs in $\cN_5$ end in $\pm1$ (see Example~\ref{ex:small}), so the last fold is not available and $\cN_{10}$ is empty. It is not difficult to keep track of the proviso in general, and one obtains the following.

\begin{corollary}[Admissible sizes]\label{cor:sizes}
The set $\cN_n$ is non-empty if and only if the binary expansion of $n$ is $1$, or begins with $11$, or consists of $10$ followed by a (possibly empty) string of ones; equivalently, if and only if $n=1$ or $3\cdot2^j-1\le n\le2^{j+2}-1$ for some $j\ge0$. In this case $|\cN_n|=2^{\lfloor\log_2n\rfloor}$.
\end{corollary}

Thus nilpotent pairs exist for $n=1,2,3$, $5\le n\le7$, $11\le n\le15$, $23\le n\le31$, and so forth, that is, for $n=1$ and for the upper half of each dyadic block $[2^{j+1},2^{j+2})$ together with the integer immediately below it, but not for $n=4$, $8\le n\le10$ or $16\le n\le22$. Since a tridiagonal matrix with non-vanishing off-diagonal entries has rank at least $n-1$, every matrix \eqref{eq:M} whose pair lies in $\cN_n$ is similar to a single Jordan block.

\begin{example}\label{ex:small}
The nilpotent pairs of size at most $7$ are $o$ and
\begin{align*}
n&=2: & d&=(-\delta,\delta), & b&=(-1);\\
n&=3: & d&=(0,0,0), & b&=(\varepsilon,-\varepsilon);\\
n&=5: & d&=(-\delta,\delta,0,\delta,-\delta), & b&=(-1,\varepsilon,-\varepsilon,-1);\\
n&=6: & d&=(0,0,-\delta,\delta,0,0), & b&=(\varepsilon,-\varepsilon,-1,-\varepsilon,\varepsilon);\\
n&=7: & d&=(0,\dots,0), & b&=(\varepsilon,-\varepsilon,\varepsilon',-\varepsilon',-\varepsilon,\varepsilon),
\end{align*}
where $\delta,\varepsilon,\varepsilon'$ range over $\{-1,1\}$. The pairs of size $2$ cannot be folded evenly, which is why $\cN_4$ is empty.
\end{example}

An even fold always creates two non-zero diagonal entries. Hence a nilpotent pair with $d=0$ is obtained from $o$ by odd folds alone, and for such pairs Theorem~\ref{thm:main} reduces to the theorem of \cite{Blattner,Shen}.

Returning to the matrices \eqref{eq:M}, we note that each pair in $\cC_n$ arises from $2^{n-1}$ of them, and that a vanishing product $a_kc_k$ arises from $5$ of the $9$ possible values of $(a_k,c_k)$. Decomposing an arbitrary matrix \eqref{eq:M} into blocks at the positions where $a_kc_k=0$, we obtain the following.

\begin{corollary}[Count]\label{cor:count}
Let $N(n)$ denote the number of nilpotent matrices of the form \eqref{eq:M}, and let $U(x):=\sum_n2^{\,n-1+\lfloor\log_2n\rfloor}x^n$, where the sum is over those $n$ for which $\cN_n$ is non-empty. Then
\[
\sum_{n\ge1}N(n)\,x^n=\frac{U(x)}{1-5\,U(x)}.
\]
In particular, $N(n)=1,9,73,585,4809,39113,318921,2599625,21186505,172688585$ for $n=1,\dots,10$.
\end{corollary}

Since a pair in $\cC_n$ is a fold precisely when its entries satisfy a few explicit equalities, Theorem~\ref{thm:main} also shows that the nilpotency of a matrix \eqref{eq:M} can be decided with $O(n)$ comparisons of its entries (Remark~\ref{rem:algo}).

\subsection{Rigidity modulo powers of two}
The proof of Theorem~\ref{thm:main} proceeds by a $2$-adic analysis of the equation $\chi_p=\lambda^n$, in the spirit of \cite{Blattner}, and in fact gives the following more precise statement, which is the main technical result of the paper. For $f\in\ZZ[\lambda]$ we let $v(f)$ denote the largest integer $j$ such that $2^j$ divides every coefficient of $f$, with the convention $v(0)=\infty$.

\begin{theorem}[Rigidity modulo $4$]\label{thm:mod4}
Let $n\ge2$, let $h:=\lfloor n/2\rfloor$, and let $P\in\cC_n$. Then $\chi_P\equiv\lambda^n\pmod 4$ if and only if $P$ is a fold of a pair $p\in\cC_h$ with $\chi_p\equiv\lambda^h\pmod 2$, that is to say $P=\Phi_\varepsilon(p)$ if $n$ is odd and $P=\Psi_\delta(p)$ if $n$ is even. In this case the pair $p$ and the sign are uniquely determined by $P$, and
\begin{equation}\label{eq:descent}
v(\chi_P-\lambda^n)=v(\chi_p-\lambda^h)+1.
\end{equation}
\end{theorem}

Thus the congruence $\chi_P\equiv\lambda^n\pmod4$ already forces $P$ to have the shape of a fold, and each fold accounts for exactly one power of two. Theorem~\ref{thm:main} follows at once (see Section~\ref{sec:proofs}), and an induction on $n$ gives the following.

\begin{corollary}[Nilpotency modulo powers of two]\label{cor:power}
Let $n\ge1$ and $\ell:=\lfloor\log_2n\rfloor$.
\begin{enumerate}
\item[(i)] If $p\in\cC_n$ and $\chi_p\equiv\lambda^n\pmod{2^{\ell+1}}$, then $\chi_p=\lambda^n$.
\item[(ii)] If $n\ge2$ and $\cN_n$ is non-empty, then there exists $p\in\cC_n$ with $v(\chi_p-\lambda^n)=\ell$. Thus the modulus in \textup{(i)} cannot be lowered.
\item[(iii)] If $\cN_n$ is empty, then the largest value of $v(\chi_p-\lambda^n)$ for $p\in\cC_n$ is $t+1$, where $t\ge0$ is the number of consecutive digits $1$ at the end of the binary expansion of $n$.
\end{enumerate}
\end{corollary}

For instance, there are pairs $p\in\cC_{19}$ with $\chi_p\equiv\lambda^{19}\pmod8$, but none with $\chi_p\equiv\lambda^{19}\pmod{16}$, and for even $n$ with $\cN_n$ empty the congruence $\chi_p\equiv\lambda^n\pmod4$ has no solutions in $\cC_n$ at all.

\begin{remark}\label{rem:full}
Corollary~\ref{cor:power} concerns pairs in $\cC_n$, that is to say matrices \eqref{eq:M} all of whose products $a_kc_k$ are non-zero. For a general matrix \eqref{eq:M}, part (i) may be applied to each of the blocks into which the matrix splits at its vanishing products, but not to the matrix as a whole. For instance, the direct sum of two matrices with pairs
\[
((-1,1),(1))\qquad\text{and}\qquad((-1,-1,0,1,1),(-1,-1,-1,-1))
\]
has characteristic polynomial $(\lambda^2-2)(\lambda^5+2\lambda^3+4\lambda)=\lambda^7-8\lambda$, which is congruent to $\lambda^7$ modulo $8=2^{1+\lfloor\log_27\rfloor}$. Let $J(n)$ denote the least integer $J$ such that every matrix $M$ of the form \eqref{eq:M} with $\det(\lambda-M)\equiv\lambda^n\pmod{2^J}$ is nilpotent. A computation gives $J(n)=1+\lfloor\log_2n\rfloor$ for $n\le6$ and for $n=8$, but $J(7)=4$ and $J(n)=5$ for $9\le n\le12$. We do not know how $J(n)$ grows.
\end{remark}

\subsection{Overview of the proof}\label{sec:overview}
We now sketch the proof of Theorem~\ref{thm:mod4}. The ``if'' direction and the identity \eqref{eq:descent} are immediate from \eqref{eq:fold}: if $q\equiv\lambda^h\pmod2$, then $q+\lambda^h$ is divisible by $2$ exactly once, so squaring $q$ improves the approximation by exactly one power of two. The content of the theorem is the ``only if'' direction, and here we follow the strategy of \cite{Blattner}, which is to study the equation $\chi_P=\lambda^n$ first modulo $2$ and then modulo $4$.

Modulo $2$ a pair in $\cC_n$ retains only the parities of its diagonal entries, and the congruence $\chi_P\equiv\lambda^n\pmod 2$ turns out to determine these parities uniquely (Lemma~\ref{lem:mod2}). As folds of nilpotent pairs are nilpotent, the solution modulo $2$ of size $n$ is then the fold of the solution of size $\lfloor n/2\rfloor$. Consequently the left half $p$ of $P$ is nilpotent modulo $2$, and its fold $Q$ satisfies $\chi_Q\equiv\lambda^n\pmod4$ and agrees with $P$ modulo $2$. (For even $n$ the last diagonal entry of $p$ has to be adjusted, and $Q$ need not have its entries in $\{-1,0,1\}$; this is the source of the gaps in Corollary~\ref{cor:sizes}.)

Modulo $4$ one can then linearise. Since $\chi$ is an affine function of each entry of the pair and $P-Q$ is divisible by $2$, the difference $\chi_P-\chi_Q$ is congruent modulo $4$ to the first-order term of its Taylor expansion around $Q$. As $Q$ is symmetric modulo $2$ under reversal, the partial derivatives of $\chi$ at mirror positions agree modulo $2$; the key point is that there are no further linear relations among them (Lemma~\ref{lem:rank}). Since $P$ and $Q$ agree on the left half, the congruence $\chi_P\equiv\chi_Q\pmod4$ therefore forces $P\equiv Q\pmod4$ in every entry.

It is only at this stage that the restriction of the diagonal entries to $\{-1,0,1\}$ is used, in the form that two such numbers which are congruent modulo $4$ are equal; the hypothesis $b_k=\pm1$ enters already in the first two steps. The examples in Remark~\ref{rem:population} show that the classification does break down for slightly larger populations.

The proofs of Theorems~\ref{thm:main} and~\ref{thm:mod4}, of Corollary~\ref{cor:sizes} and of Corollary~\ref{cor:power}(i) have been formalised in Lean~4; see Section~\ref{sec:lean}.

\subsection{Relation to previous work}
As the above sketch indicates, the tools used in the proof are not new. The identities \eqref{eq:fold} are a form of the folding lemma of Mend\`es France \cite{MendesFrance} (see also \cite{vdPShallit}), under the familiar dictionary between tridiagonal matrices and continued fractions. The uniqueness statement modulo $2$ is, in different language, part of the theory of binary sequences with a perfect linear complexity profile \cite{AHN,BaumSweet,Blackburn,WangMassey}. The linear independence used modulo $4$ is a statement of the kind for which Garnett and Shader introduced the nilpotent-centralizer method \cite{GarnettShader1,GarnettShader2} in the study of spectrally arbitrary sign patterns; in that literature one finds nilpotent tridiagonal matrices with prescribed sign patterns and with real entries of unrestricted magnitude \cite{Behn,Drew}, whereas here all magnitudes are equal to $1$. The contribution of the present paper is the observation that these tools combine to give a complete description of the sets $\cN_n$. We are not aware of an earlier appearance of Theorem~\ref{thm:main} or of its corollaries; in particular, the sequence in Corollary~\ref{cor:count} was not in the OEIS at the time of writing, the analogous sequence for unrestricted matrices with entries in $\{-1,0,1\}$ being \cite{OEIS}.

Matrices with pairs in $\cN_n$ have nevertheless appeared in the literature in two other contexts: as the test matrices of Liu for nonsymmetric tridiagonal eigensolvers, and in a recent design of exceptional points of high order in photonics \cite{Takata}. We discuss these connections in Section~\ref{sec:liu}.

\subsection{The eigenvalue zero in Corless's family}
In Section~\ref{sec:mult} we return to the skew-symmetric matrices of Corless, and to a question left open in \cite[Section~7]{Blattner}: which multiplicities of the eigenvalue zero occur in dimension $m$? The $2$-adic method shows that the multiplicity is at most $2^{v_2(m+1)}-1$, where $v_2$ denotes the $2$-adic valuation, and that this bound is attained for every $m$ (Proposition~\ref{prop:bound}). When $m=2^K-1$ it also identifies the matrices with exactly two non-zero eigenvalues, counted with algebraic multiplicity, which arise from the nilpotent ones by an ``anti-fold'' (Proposition~\ref{prop:deg1}). Computations for $m\le63$ suggest that for $K\ge4$ the multiplicity $m-6$ is absent, although $m-4$ occurs (Conjecture~\ref{conj:deg3}).

\section{Preliminaries}\label{sec:prelim}

\subsection{Continuants}\label{sec:continuants}
Let $p=(d,b)$ be a pair of size $n$ over a commutative ring. For $1\le i\le j\le n$ let $\seg ij$ denote the characteristic polynomial of the pair formed by the entries $d_i,\dots,d_j$ and $b_i,\dots,b_{j-1}$, with the conventions $\seg{i}{i-1}:=1$ and $\seg{i}{i-2}:=0$. We write
\[
\ell_k:=\seg1k\qquad\text{and}\qquad r_k:=\seg kn
\]
for the characteristic polynomials of the first $k$ and of the last $n+1-k$ positions; thus $\chi_p=\ell_n=r_1$, while $\ell_0=r_{n+1}=1$ and $r_{n+2}=0$. Laplace expansion gives
\begin{align}
\chi_p&=\ell_k\,r_{k+1}-b_k\,\ell_{k-1}\,r_{k+2},\label{eq:edge}\\
\chi_p&=(\lambda-d_k)\,\ell_{k-1}\,r_{k+1}-b_{k-1}\,\ell_{k-2}\,r_{k+1}-b_k\,\ell_{k-1}\,r_{k+2}\label{eq:vertex}
\end{align}
for $1\le k\le n-1$ and $1\le k\le n$ respectively, where terms involving $b_0$ or $b_n$ are to be omitted. Applied to the last $n+1-k$ positions, the case $k=1$ of \eqref{eq:vertex} is the three-term recurrence
\begin{equation}\label{eq:three}
r_k=(\lambda-d_k)\,r_{k+1}-b_k\,r_{k+2}\qquad(1\le k\le n).
\end{equation}
It follows from \eqref{eq:edge} and \eqref{eq:vertex} that $\chi_p$ is unchanged when $p$ is reversed, and that it is an affine function of each $d_k$ and of each $b_k$, with
\begin{equation}\label{eq:derivatives}
\frac{\partial\chi_p}{\partial d_k}=-A_k,\qquad\frac{\partial\chi_p}{\partial b_k}=-B_k,
\end{equation}
where
\[
A_k:=\ell_{k-1}\,r_{k+1},\qquad B_k:=\ell_{k-1}\,r_{k+2}.
\]
We also record the identity
\begin{equation}\label{eq:continuant}
r_2\,\ell_{n-1}-\chi_p\,\seg2{n-1}=b_1b_2\cdots b_{n-1},
\end{equation}
which follows by induction on $n$ from the case $k=n$ of \eqref{eq:vertex}.

\subsection{Folding}
We now verify \eqref{eq:fold}.

\begin{lemma}[Folding]\label{lem:fold}
Let $p$ be a pair of size $n$ over a commutative ring. Then $\chi_{\Phi_\varepsilon(p)}=\lambda\chi_p^2$ and $\chi_{\Psi_\delta(p)}=\chi_p^2$ for all $\varepsilon,\delta\in\{-1,1\}$.
\end{lemma}

\begin{proof}
Write $\chi:=\chi_p$ and $\chi^-:=\seg1{n-1}$. For the odd fold we apply \eqref{eq:vertex} at the central position $n+1$. By reversal invariance, the two halves of $\Phi_\varepsilon(p)$ have characteristic polynomial $\chi$, and after deletion of the position adjacent to the centre they have characteristic polynomial $\chi^-$; hence $\chi_{\Phi_\varepsilon(p)}=\lambda\chi^2-\varepsilon\chi^-\chi+\varepsilon\chi\chi^-=\lambda\chi^2$. For the even fold we apply \eqref{eq:edge} at the central edge. By \eqref{eq:derivatives}, the two halves of $\Psi_\delta(p)$ have characteristic polynomials $\chi+\delta\chi^-$ and $\chi-\delta\chi^-$, and after deletion of the two central positions both have characteristic polynomial $\chi^-$; hence $\chi_{\Psi_\delta(p)}=(\chi+\delta\chi^-)(\chi-\delta\chi^-)+(\chi^-)^2=\chi^2$.
\end{proof}

\begin{remark}\label{rem:cf}
If $M$ is a tridiagonal matrix with pair $p$, then $e_1^\tp(\lambda-M)^{-1}e_1=\seg2n/\seg1n$ is the finite continued fraction with partial quotients $\lambda-d_1,\dots,\lambda-d_n$ and partial numerators $-b_1,\dots,-b_{n-1}$. In this language Lemma~\ref{lem:fold} is the folding lemma \cite{MendesFrance,vdPShallit}: when the $b_k$ are $\pm1$, passing from $p$ to $\Phi_\varepsilon(p)$ changes this quantity by $\pm1/(\lambda\chi_p^2)$, and passing to $\Psi_\delta(p)$ changes it by $\pm1/\chi_p^2$. We return to this point in Proposition~\ref{prop:resolvent}.
\end{remark}

\subsection{The situation modulo \texorpdfstring{$2$}{2}}
Modulo $2$ all the products $b_k=\pm1$ become equal to $1$, and we are led to consider the pairs $(d,\one)$ over $\FF_2$, where $\one:=(1,\dots,1)$; these are the pairs of the symmetric tridiagonal matrices over $\FF_2$ with unit off-diagonal entries.

\begin{lemma}[Uniqueness modulo $2$]\label{lem:mod2}
For every $n\ge1$ there is exactly one $\nu_n\in\FF_2^n$ for which the pair $(\nu_n,\one)$ over $\FF_2$ is nilpotent. We have $\nu_1=(0)$, and for every $h\ge1$,
\begin{equation}\label{eq:nu}
\nu_{2h+1}=(\nu_h,\,0,\,\rev\nu_h)\quad\text{and}\quad\nu_{2h}=(e_1,\dots,e_{h-1},\,e_h+1,\,e_h+1,\,e_{h-1},\dots,e_1),
\end{equation}
where $e:=\nu_h$. In particular $\nu_n$ is a palindrome, and its first entry is $1$ if and only if the binary expansion of $n$ begins with $10$.
\end{lemma}

In the language of continued fractions, the lemma asserts that there is exactly one rational function over $\FF_2$ with denominator $\lambda^n$ all of whose partial quotients are linear. In that form it is known: the identity \eqref{eq:riccati} below is a finite form of a criterion of Baum and Sweet \cite{BaumSweet}, and the relations \eqref{eq:moments} are a finite form of those of Wang and Massey \cite{WangMassey}; see the introduction of \cite{AHN} for the history, and \cite{Blackburn} for a closely related count. As the argument is short, we include it.

\begin{proof}
By Lemma~\ref{lem:fold}, applied over $\FF_2$, if a pair $(e,\one)$ of size $h$ is nilpotent, then so are the pairs of sizes $2h+1$ and $2h$ whose diagonals are obtained from $e$ as in \eqref{eq:nu}. This proves the existence of $\nu_n$, and it shows that \eqref{eq:nu} follows from the uniqueness. The two final assertions follow from \eqref{eq:nu}: both vectors there are palindromes, and for $n\ge4$ the first entry of $\nu_n$ is that of $\nu_{\lfloor n/2\rfloor}$, while $\nu_2=(1,1)$ and $\nu_3=(0,0,0)$.

It remains to prove the uniqueness. Let $d\in\FF_2^n$ be arbitrary, and let $r_1,\dots,r_{n+2}$ be the polynomials of Section~\ref{sec:continuants} for the pair $(d,\one)$, so that $r_k=(\lambda+d_k)r_{k+1}+r_{k+2}$ by \eqref{eq:three}. We claim that
\begin{equation}\label{eq:riccati}
r_k^2+\lambda\,r_kr_{k+1}+r_{k+1}^2+1=(1+\lambda)\,g_k^2,\qquad\text{where }g_k:=\sum_{j=k}^{n}d_j\,r_{j+1},
\end{equation}
for $1\le k\le n+1$. For $k=n+1$ both sides vanish. If \eqref{eq:riccati} holds for $k+1$, we substitute the recurrence for $r_k$ and find that the left-hand side for $k$ exceeds the left-hand side for $k+1$ by $d_k(1+\lambda)r_{k+1}^2$; since $d_k^2=d_k$ and $g_k=g_{k+1}+d_kr_{k+1}$, this gives \eqref{eq:riccati} for $k$.

Now suppose that $(d,\one)$ is nilpotent, so that $r_1=\lambda^n$, and write $r_2=\sum_{j=1}^{n}c_j\lambda^{n-j}$, with $c_j:=0$ for $j>n$. We apply \eqref{eq:riccati} with $k=1$. In a polynomial of the form $(1+\lambda)g^2$ over $\FF_2$, the coefficients of $\lambda^{2i}$ and of $\lambda^{2i+1}$ agree for every $i$. Comparing the coefficients of $\lambda^{2(n-i)}$ and of $\lambda^{2(n-i)+1}$ in $r_2^2+\lambda^{n+1}r_2+\lambda^{2n}+1$, we obtain
\begin{equation}\label{eq:moments}
c_n=1\qquad\text{and}\qquad c_i=c_{2i}+c_{2i+1}\quad\text{for }1\le i<n.
\end{equation}
These relations determine $c_n,c_{n-1},\dots,c_1$ in turn, so that $r_2$ does not depend on $d$. On the other hand $r_1$ and $r_2$ determine $d$: by \eqref{eq:three}, the polynomials $\lambda+d_k$ and $r_{k+2}$ are the quotient and the remainder in the division of $r_k$ by $r_{k+1}$, for $k=1,\dots,n$ in turn.
\end{proof}

\subsection{Transversality}
Our final tool asserts that the characteristic polynomial map $p\mapsto\chi_p$ is a submersion wherever the $b_k$ do not vanish, and identifies the linear relations among its partial derivatives \eqref{eq:derivatives} when $p$ is a palindrome.

\begin{lemma}[Transversality]\label{lem:rank}
Let $p=(d,b)$ be a pair of size $n$ over a field $F$, with all $b_k$ non-zero.
\begin{enumerate}
\item[(i)] The polynomials $A_1,\dots,A_n,B_1,\dots,B_{n-1}$ span the space of polynomials in $F[\lambda]$ of degree less than $n$.
\item[(ii)] Suppose in addition that $d=\rev d$ and $b=\rev b$. Then $A_k=A_{n+1-k}$ and $B_k=B_{n-k}$ for all $k$, and if $t\in F^n$ and $u\in F^{n-1}$ are such that $\sum_kt_kA_k+\sum_ku_kB_k=0$, then $t+\rev t=0$ and $u+\rev u=0$. For odd $n$ the central entry of $t$ vanishes, and for even $n$ the central entry of $u$ vanishes.
\end{enumerate}
\end{lemma}

\begin{proof}
Let $W$ be the span of the polynomials in (i). We show by induction on $k$ that $W$ contains $f\,r_{k+1}$ for every polynomial $f$ of degree less than $k$, for $1\le k\le n$. As $r_{n+1}=1$, the case $k=n$ is the assertion of (i). For $k=1$ the claim holds because $r_2=A_1$. Now let $1\le k<n$, and suppose that the claim holds up to $k$. Since $\ell_k$ and $\ell_{k-1}$ are monic of degrees $k$ and $k-1$, every polynomial of degree at most $k$ is of the form $\alpha\ell_k+\beta\ell_{k-1}+f$ with $\alpha,\beta\in F$ and $\deg f<k-1$. Here $\ell_kr_{k+2}=A_{k+1}$ and $\ell_{k-1}r_{k+2}=B_k$ lie in $W$, and by \eqref{eq:three},
\[
b_k\,f\,r_{k+2}=(\lambda-d_k)f\,r_{k+1}-f\,r_k,
\]
where both terms on the right lie in $W$ by the induction hypothesis (for $k=1$ we have $f=0$). As $b_k\ne0$, this proves the claim for $k+1$.

In the situation of (ii), reversal invariance gives $A_k=A_{n+1-k}$ and $B_k=B_{n-k}$. By (i), the $n$ polynomials $A_k$ with $k\le\lceil n/2\rceil$ and $B_k$ with $k\le\lfloor n/2\rfloor$ span a space of dimension $n$, and hence are linearly independent. Collecting terms in the given relation, we obtain (ii).
\end{proof}

Part (i) says that the Jacobian of the map $p\mapsto\chi_p$ has full rank. Statements of this kind are the object of the nilpotent-centralizer method of Garnett and Shader \cite{GarnettShader1,GarnettShader2}, who applied it to a tridiagonal sign pattern with non-zero diagonal entries at the two ends only; for tridiagonal matrices the three-term recurrence permits the direct argument above. When $d=0$ and $b$ is constant, the $B_k$ correspond to the polynomials $c_e$ of \cite{Blattner}, and part (ii) recovers \cite[Lemma~3.1]{Blattner}, in every dimension and over every field.

\section{Proofs of the main results}\label{sec:proofs}

\subsection{Proof of Theorem~\ref{thm:mod4}}
We begin with a computation. Let $p$ be a pair of size $h$ with integer entries such that $q:=\chi_p\equiv\lambda^h\pmod2$, and let $P$ be a fold of $p$, of size $n$. By Lemma~\ref{lem:fold},
\[
\chi_P-\lambda^n=\lambda^{n-2h}\,(q-\lambda^h)(q+\lambda^h).
\]
Since the coefficients of $q-\lambda^h$ are even, $(q+\lambda^h)/2=\lambda^h+(q-\lambda^h)/2$ is a monic polynomial with integer coefficients, and so $v(q+\lambda^h)=1$. As $v$ is additive on products by Gauss's lemma, we obtain \eqref{eq:descent}, and in particular $\chi_P\equiv\lambda^n\pmod4$. This proves the ``if'' direction of the theorem. The uniqueness of $p$ and of the sign is clear, as both can be read off from a fold $P=(d,b)\in\cC_n$ of a pair $p\in\cC_h$:
\begin{align}
p&=((d_1,\dots,d_h),\,(b_1,\dots,b_{h-1})),&\varepsilon&=b_h&&\text{if $n$ is odd},\label{eq:parentodd}\\
p&=((d_1,\dots,d_{h-1},0),\,(b_1,\dots,b_{h-1})),&\delta&=d_{h+1}&&\text{if $n$ is even}.\label{eq:parenteven}
\end{align}

Now suppose conversely that $P=(d,b)\in\cC_n$ satisfies $\chi_P\equiv\lambda^n\pmod4$.

\emph{Step 1: reduction modulo $2$.} The reduction of $P$ modulo $2$ is nilpotent, so $d\equiv\nu_n\pmod 2$ by Lemma~\ref{lem:mod2}.

\emph{Step 2: the fold of the left half.} We define an auxiliary pair $p$ of size $h$ with integer entries, which for even $n$ need not lie in $\cC_h$, and a fold $Q$ of $p$, as follows. If $n=2h+1$, let
\[
p:=((d_1,\dots,d_h),\,(b_1,\dots,b_{h-1}))\qquad\text{and}\qquad Q:=\Phi_\varepsilon(p)\quad\text{with }\varepsilon:=b_h.
\]
If $n=2h$, let
\[
p:=((d_1,\dots,d_{h-1},\,d_h+1),\,(b_1,\dots,b_{h-1}))\qquad\text{and}\qquad Q:=\Psi_1(p),
\]
so that
\[
Q=\bigl((d_1,\dots,d_h,\,d_h+2,\,d_{h-1},\dots,d_1),\ (b_1,\dots,b_{h-1},\,-1,\,b_{h-1},\dots,b_1)\bigr).
\]
In both cases Step 1 and the recursion \eqref{eq:nu} show that the diagonal of $p$ is congruent to $\nu_h$ modulo $2$, and that the diagonal of $Q$ is congruent to $\nu_n$, and hence to $d$. Thus $\chi_p\equiv\lambda^h\pmod2$, and therefore $\chi_Q\equiv\lambda^n\pmod4$ by the computation above. By construction, $Q$ agrees with $P$ in the entries $d_1,\dots,d_h$ and $b_1,\dots,b_{h-1}$, and also in $b_h$ if $n$ is odd.

\emph{Step 3: linearisation modulo $4$.} Write $Q=(d',b')$. By Step 2 the differences $d_k-d'_k$ and $b_k-b'_k$ are all even. Since $\chi$ is affine in each of its arguments, \eqref{eq:derivatives} gives
\[
\chi_P\equiv\chi_Q-\sum_{k=1}^{n}(d_k-d'_k)A_k-\sum_{k=1}^{n-1}(b_k-b'_k)B_k\pmod 4,
\]
where $A_k$ and $B_k$ are evaluated at $Q$. Both characteristic polynomials are congruent to $\lambda^n$ modulo $4$. Dividing by $2$ and reducing modulo $2$, we obtain a relation $\sum_kt_k\bar A_k+\sum_ku_k\bar B_k=0$ over $\FF_2$, in which $t_k$ and $u_k$ are the residues of $(d_k-d'_k)/2$ and $(b_k-b'_k)/2$, and $\bar A_k$, $\bar B_k$ are the polynomials $A_k$, $B_k$ of the pair $(\nu_n,\one)$, the reduction of $Q$ modulo $2$. This pair is a palindrome, so Lemma~\ref{lem:rank}(ii) applies: we have $t_k=t_{n+1-k}$ and $u_k=u_{n-k}$ for all $k$, the central entry of $t$ vanishes if $n$ is odd, and the central entry of $u$ vanishes if $n$ is even. But $t_k=0$ for $k\le h$ by Step 2, and likewise $u_k=0$ for $k<h$, and for $k=h$ if $n$ is odd. Hence $t=0$ and $u=0$, that is to say
\begin{equation}\label{eq:PQ}
d_k\equiv d'_k\qquad\text{and}\qquad b_k\equiv b'_k\pmod4\qquad\text{for all }k.
\end{equation}

\emph{Step 4: the shape of $P$.} Two elements of $\{-1,0,1\}$ which are congruent modulo $4$ are equal. If $n$ is odd, then $Q\in\cC_n$, so that \eqref{eq:PQ} gives $P=Q=\Phi_\varepsilon(p)$, as required. If $n$ is even, then all entries of $Q$ other than $d'_{h+1}=d_h+2$ lie in $\{-1,0,1\}$, and \eqref{eq:PQ} gives $d_k=d_{n+1-k}$ and $b_k=b_{n-k}$ for $k<h$, as well as $b_h=-1$ and $d_{h+1}\equiv d_h+2\pmod4$. The last congruence forces $d_{h+1}=-d_h\ne0$. Hence $P=\Psi_\delta(p')$, where $\delta:=d_{h+1}$ and $p':=((d_1,\dots,d_{h-1},0),(b_1,\dots,b_{h-1}))\in\cC_h$. Finally $p'\equiv p\pmod2$, so that $\chi_{p'}\equiv\lambda^h\pmod2$.\qed

\subsection{Proof of Theorem~\ref{thm:main} and its corollaries}
\begin{proof}[Proof of Theorem~\ref{thm:main}]
The sets on the right-hand sides are contained in $\cN_{2h+1}$ and $\cN_{2h}$ by Lemma~\ref{lem:fold} and the remarks following \eqref{eq:fold}. Conversely, let $P\in\cN_n$ with $n\ge2$. By Theorem~\ref{thm:mod4}, $P$ is a fold of a unique pair $p\in\cC_h$, with a unique sign, and $\lambda^{n-2h}\chi_p^2=\lambda^n$ forces $\chi_p=\lambda^h$. If $n$ is even, the last diagonal entry of $p$ vanishes, since otherwise $\Psi_\delta(p)$ would not lie in $\cC_n$.
\end{proof}

\begin{remark}[Recognition]\label{rem:algo}
Theorem~\ref{thm:main} yields a test for nilpotency which avoids the computation of the characteristic polynomial. A pair $P=(d,b)\in\cC_n$ with $n\ge2$ is a fold of a pair in $\cC_h$, where $h:=\lfloor n/2\rfloor$, if and only if $d_k=d_{n+1-k}$ and $b_k=b_{n-k}$ for all $k<h$, and in addition
\[
d_h=d_{h+2},\quad d_{h+1}=0,\quad b_{h+1}=-b_h\qquad\text{or}\qquad d_{h+1}=-d_h\ne0,\quad b_h=-1,
\]
according as $n$ is odd or even. If these conditions fail, $P$ is not nilpotent. If they hold, then $P$ is nilpotent if and only if the pair $p$ in \eqref{eq:parentodd} or \eqref{eq:parenteven} is, by \eqref{eq:fold}, and one repeats the test with $p$ in place of $P$. As the sizes at least halve at each step, this decides the nilpotency of $P$ with $O(n)$ comparisons of entries, and after splitting at the vanishing products the same holds for an arbitrary matrix \eqref{eq:M}. When the test succeeds it also returns the signs of the successive folds, that is, the description of $P$ furnished by Theorem~\ref{thm:main}.
\end{remark}

\begin{proof}[Proof of Corollary~\ref{cor:sizes}]
By Lemma~\ref{lem:mod2}, the last diagonal entry of a pair in $\cN_h$ is odd if the binary expansion of $h$ begins with $10$, and is even, hence zero, otherwise. Theorem~\ref{thm:main} therefore gives $|\cN_{2h+1}|=2|\cN_h|$ for all $h\ge1$, while $|\cN_{2h}|=2|\cN_h|$ if the expansion of $h$ does not begin with $10$, and $\cN_{2h}=\emptyset$ if it does. Thus, along the sequence $1=\lfloor n/2^\ell\rfloor,\dots,\lfloor n/2\rfloor,n$ of binary prefixes of $n$, the cardinality doubles at each step, unless at some step a digit $0$ is appended to a prefix beginning with $10$, in which case $\cN_n$ is empty. The latter happens precisely when the expansion of $n$ begins with $10$ and contains a further digit $0$.
\end{proof}

\begin{proof}[Proof of Corollary~\ref{cor:count}]
Let $n_1,\dots,n_j$ be the sizes of the blocks into which a matrix \eqref{eq:M} is cut at the positions with $a_kc_k=0$. By \eqref{eq:edge} the matrix is nilpotent if and only if each block is. The number of nilpotent matrices with given block sizes is therefore $5^{j-1}\prod_i2^{n_i-1}|\cN_{n_i}|$, and summing over all compositions of $n$ we obtain $\sum_nN(n)x^n=\sum_{j\ge1}5^{j-1}U(x)^j$.
\end{proof}

\subsection{Proof of Corollary~\ref{cor:power}}
We prove (i) by induction on $n$, the case $n=1$ being trivial. Let $n\ge2$ and suppose that $v(\chi_P-\lambda^n)\ge\ell+1$ for some $P\in\cC_n$. As $\ell\ge1$, Theorem~\ref{thm:mod4} applies and produces a pair $p\in\cC_h$ with $v(\chi_p-\lambda^h)\ge\ell=\lfloor\log_2h\rfloor+1$. By the induction hypothesis $\chi_p=\lambda^h$, and hence $\chi_P=\lambda^n$.

For (ii), we note that the sequence of binary prefixes of $n$ passes through $2$ or through $3$. The pairs $((-1,1),(1))$ and $((0,0,0),(1,1))$ have characteristic polynomials $\lambda^2-2$ and $\lambda^3-2\lambda$. The first has the same diagonal as one of the two elements of $\cN_2$ (the other has the reversed diagonal), and the second has the same diagonal as both elements of $\cN_3$. As the diagonal of a fold depends only on the diagonal of its parent and on the sign, the folds prescribed by the remaining $\ell-1$ digits of $n$ may be applied to these pairs without leaving the classes $\cC_k$, and by \eqref{eq:descent} each of these folds increases $v(\chi-\lambda^k)$ by exactly one. The resulting pair $p\in\cC_n$ has $v(\chi_p-\lambda^n)=\ell$.

For (iii), write $V(m)$ for the largest value of $v(\chi_p-\lambda^m)$ over $p\in\cC_m$, so that $V(m)=\infty$ precisely when $\cN_m$ is non-empty. By Lemma~\ref{lem:mod2}, the pairs $p\in\cC_m$ with $v(\chi_p-\lambda^m)\ge1$ are those whose diagonal is congruent to $\nu_m$ modulo $2$. In particular $V(m)\ge1$, and the last diagonal entry of such a pair is $\pm1$ if the binary expansion of $m$ begins with $10$, and is $0$ otherwise. Now let $m\ge2$ and $h:=\lfloor m/2\rfloor$. By Theorem~\ref{thm:mod4}, the pairs $P\in\cC_m$ with $v(\chi_P-\lambda^m)\ge2$ are the folds, lying in $\cC_m$, of the pairs $p\in\cC_h$ with $v(\chi_p-\lambda^h)\ge1$, and \eqref{eq:descent} holds for them. If $m$ is odd, or if $m$ is even and the expansion of $h$ does not begin with $10$, every such $p$ has a fold in $\cC_m$, and so $V(m)=V(h)+1$. If $m$ is even and the expansion of $h$ begins with $10$, no such $p$ has a fold in $\cC_m$, and $V(m)=1$. Finally, if $\cN_n$ is empty, then the binary expansion of $n$ begins with $10$ and contains a further digit $0$. Deleting the $t$ final digits $1$ leaves an even number $2h$ such that the expansion of $h$ begins with $10$, and we conclude that $V(n)=t+V(2h)=t+1$.\qed

\subsection{The resolvent, and the role of the population}
For a nilpotent pair the continued fraction of Remark~\ref{rem:cf} can be evaluated in closed form. For $d=0$ and $n=2^K-1$ the following proposition is \cite[Corollary~4.3]{Blattner}.

\begin{proposition}[Resolvent]\label{prop:resolvent}
Let $n$ be such that $\cN_n$ is non-empty, and let $\ell:=\lfloor\log_2n\rfloor$. For every $p\in\cN_n$ and every tridiagonal matrix $M$ with pair $p$,
\[
e_1^\tp(\lambda-M)^{-1}e_1=\sum_{t=0}^{\ell}\eta_t\,\lambda^{-\lfloor n/2^t\rfloor}
\]
for some signs $\eta_t=\eta_t(p)\in\{-1,1\}$ with $\eta_\ell=1$. Moreover, the map $p\mapsto(\eta_0,\dots,\eta_{\ell-1})$ is a bijection from $\cN_n$ to $\{-1,1\}^\ell$.
\end{proposition}

\begin{proof}
For a pair $p=(d,b)$ of size $n$ put $w(p):=r_2/\chi_p$ and $\pi:=b_1\cdots b_{n-1}$. Arguing as in the proof of Lemma~\ref{lem:fold}, with the first position removed and with \eqref{eq:continuant} in place of the cancellation there, one finds
\[
w(\Phi_\varepsilon(p))=w(p)+\frac{\varepsilon\pi}{\lambda\chi_p^2},\qquad w(\Psi_\delta(p))=w(p)-\frac{\delta\pi}{\chi_p^2}.
\]
For $p\in\cN_n$ the two increments are $\pm\lambda^{-(2n+1)}$ and $\pm\lambda^{-2n}$, and each sign occurs for exactly one value of $\varepsilon$ or of $\delta$. The claims now follow from Theorem~\ref{thm:main} by induction on $n$, starting from $w(o)=\lambda^{-1}$.
\end{proof}

\begin{remark}[The role of the population]\label{rem:population}
Iterated folds of $o$ with integer entries exist in every dimension; for instance $\Psi_1(\Psi_1(o))$ has $d=(-1,0,2,-1)$ and $b=(-1,-1,-1)$. Theorem~\ref{thm:main} asserts that, as long as the diagonal entries are confined to $\{-1,0,1\}$ and the products to $\{-1,1\}$, these are the only nilpotent pairs. Neither restriction can be relaxed very far. The pair with
\[
d=(0,-2,1,-1,2,0),\qquad b=(-1,-1,-1,-1,-1)
\]
is nilpotent and is not a fold, and with $\omega:=e^{2\pi i/3}$ the pair with
\[
d=(-1,1,-1,1),\qquad b=(\omega,-1,\omega^2)
\]
has characteristic polynomial $(\lambda^2-1)^2+2(\lambda^2-1)+1=\lambda^4$, although all its products have modulus $1$ and $\cN_4$ is empty.
\end{remark}

\section{Two appearances of folds in the literature}\label{sec:liu}

\subsection{The test matrices of Liu}
Liu \cite{Liu} devised an algorithm which produces, in every dimension, tridiagonal matrices with non-vanishing off-diagonal entries whose only eigenvalue is $0$ (see \cite{FerreiraParlett}). These ``one-point spectrum'' matrices consist of a single Jordan block, and have become standard test cases for nonsymmetric tridiagonal eigensolvers \cite{Berhanu,BGT,FerreiraParlett,FerreiraParlettDopico}. Three of the instances in the literature have entries in $\{-1,0,1\}$: one of size $6$ in \cite{FerreiraParlett}, one of size $14$ in \cite{FerreiraParlettDopico}, and one of size $28$ in \cite[Section~5]{BGT}. By Theorem~\ref{thm:main} their pairs, which we denote by $\mathrm{Liu}_n$, must be iterated folds of $o$, and indeed one finds
\[
\mathrm{Liu}_6=\Psi_1\Phi_{-1}(o),\qquad\mathrm{Liu}_{14}=\Psi_1\Phi_1\Phi_{-1}(o),\qquad\mathrm{Liu}_{28}=\Psi_1(\mathrm{Liu}_{14});
\]
thus for example $\mathrm{Liu}_6$ has $d=(0,0,-1,1,0,0)$ and $b=(-1,1,-1,1,-1)$. Not all of Liu's matrices are of this type: the instance of size $5$ recorded in \cite[Section~7.6]{Berhanu} has $d=0$ and $b=(-1,2,-\tfrac12,-\tfrac12)$.

Conversely, Theorem~\ref{thm:main} and Corollary~\ref{cor:sizes} describe all test matrices of this kind with entries in $\{-1,0,1\}$: they exist precisely in the dimensions of Corollary~\ref{cor:sizes}, where there are $2^{\lfloor\log_2n\rfloor}$ of them up to conjugation by diagonal matrices with entries $\pm1$, and Proposition~\ref{prop:resolvent} supplies an entry of the resolvent of each in closed form.

\subsection{Exceptional points}\label{sec:ep}
In the physics literature, a matrix whose Jordan normal form has a block of size $n$ is said to have an \emph{exceptional point of order $n$} at the corresponding eigenvalue. Thus an $n\times n$ matrix has an exceptional point of order $n$ at the eigenvalue zero precisely when it is nilpotent and consists of a single Jordan block. Takata, Mock, Notomi and Shinya \cite{Takata} recently proposed an inductive design of such matrices: starting from a complex symmetric $n\times n$ matrix $H$ with an exceptional point of order $n$ at the eigenvalue zero, they form the matrix in part (ii) of the proposition below and show, under an additional hypothesis on $H$, that it has an exceptional point of order $2n$ at the eigenvalue zero \cite[Theorem~1]{Takata}. For tridiagonal $H$ this construction is the even fold \eqref{eq:even}, with the sign $\delta$ replaced by the scalar $-\alpha$ and the central product $-1$ by $\beta^2=-\alpha^2$. The odd fold provides a companion to it, and both rest on the following identities, which are valid for arbitrary square matrices.

\begin{proposition}\label{prop:ep}
Let $H$ be an $n\times n$ matrix over a commutative ring, with characteristic polynomial $\chi$, and let $R$ be the $n\times n$ exchange matrix.
\begin{enumerate}
\item[(i)] If $\beta^2+\beta'^2=0$, then the matrix
\[
\begin{pmatrix} H & \beta e_n & 0\\ \beta e_n^\tp & 0 & \beta'e_1^\tp\\ 0 & \beta'e_1 & RHR\end{pmatrix}
\]
of size $2n+1$ has characteristic polynomial $\lambda\chi^2$.
\item[(ii)] If $\alpha^2+\beta^2=0$, then the matrix
\[
\begin{pmatrix} H+\alpha e_ne_n^\tp & \beta e_ne_1^\tp\\ \beta e_1e_n^\tp & RHR-\alpha e_1e_1^\tp\end{pmatrix}
\]
of size $2n$ has characteristic polynomial $\chi^2$.
\end{enumerate}
\end{proposition}

\begin{proof}
Let $\chi^-$ denote the characteristic polynomial of $H$ with its last row and column removed, which is also that of $RHR$ with its first row and column removed. Expanding along the central row and column in (i), one obtains $\lambda\chi^2-(\beta^2+\beta'^2)\chi\chi^-$. In (ii) the two diagonal blocks have characteristic polynomials $\chi-\alpha\chi^-$ and $\chi+\alpha\chi^-$, and the coupling contributes $-\beta^2(\chi^-)^2$, for a total of $\chi^2-(\alpha^2+\beta^2)(\chi^-)^2$.
\end{proof}

When $H$ is a complex tridiagonal matrix with non-vanishing off-diagonal entries and $\beta\ne0$, the matrices in (i) and (ii) are again of this form, so that if $H$ is nilpotent they are automatically single Jordan blocks. Starting from the $1\times1$ zero matrix and following the binary digits of $n$, one thus obtains for every $n$ a complex symmetric tridiagonal $n\times n$ matrix with an exceptional point of order $n$ at the eigenvalue zero, depending on $\lfloor\log_2n\rfloor$ free non-zero parameters. Amusingly, one of the designs of \cite{Takata} has already appeared above: the six-site matrix $H_{6P}$ of \cite[Eq.~(12)]{Takata} has, for the values of its parameters singled out there, an exceptional point of order $6$ at an eigenvalue $i\gamma$, and the matrix $(H_{6P}-i\gamma I)/\gamma$ then has $d=(0,0,1,-1,0,0)$ and $b=(-1,1,-1,1,-1)$. It is therefore diagonally similar to Liu's matrix of size $6$ with its rows and columns in reverse order.

Theorem~\ref{thm:main} may be read as a constraint on designs of this kind. Fix $g\ne0$, and consider the tridiagonal matrices whose diagonal entries lie in $\{0,g,-g\}$ and whose products $a_kc_k$ lie in $\{g^2,-g^2\}$. Such a matrix has an exceptional point of order $n$ at the eigenvalue zero only if $n$ is as in Corollary~\ref{cor:sizes}, and then only if its pair, after division of the matrix by $g$, is one of the $2^{\lfloor\log_2n\rfloor}$ elements of $\cN_n$. To reach the remaining orders, such as $4$, $8$, $9$, $10$ or $16,\dots,22$, one must introduce a second magnitude on the diagonal or off it, or products with other phases as in Remark~\ref{rem:population}.

\section{The eigenvalue zero in Corless's family}\label{sec:mult}

We adopt the notation of \cite{Blattner}. For $u\in\{1,i\}^{m-1}$ let $T_m(u)$ denote the $m\times m$ tridiagonal matrix with zero diagonal, superdiagonal $u$ and subdiagonal $-u$, and put $s_k:=u_k^2\in\{-1,1\}$. Then $\det(\lambda-T_m(u))=\lambda^mD_s(\lambda^{-2})$, where $D_s\in\ZZ[x]$ is the matching polynomial of the path on $m$ vertices with edge weights $s_k$; explicitly, $D_s=D_m$, where $D_0=D_1=1$ and $D_k=D_{k-1}+s_{k-1}xD_{k-2}$ \cite[Lemma~2.1]{Blattner}. In particular the eigenvalue zero of $T_m(u)$ has multiplicity
\[
\mu(s):=m-2\deg D_s.
\]
In this notation, the main result of \cite{Blattner} states that $D_s=1$ precisely when $m=2^K-1$ for some $K\ge1$ and $s$ belongs to the set $F_K$ of the $2^{K-1}$ sequences obtained from the empty sequence by $K-1$ applications of the rule $\sigma\mapsto(\sigma,\varepsilon,-\varepsilon,\rev\sigma)$ \cite[Theorem~2.2]{Blattner}. Corless observed that for $m=15$ the multiplicity $9$ does not occur \cite[Section~5]{Corless}. We begin with a sharp upper bound.

\begin{proposition}\label{prop:bound}
For every $m\ge1$ and every $s\in\{-1,1\}^{m-1}$ we have
\[
\mu(s)\le2^{v_2(m+1)}-1,
\]
where $v_2$ denotes the $2$-adic valuation, and for every $m$ equality holds for some $s$. In particular, $0$ is a simple eigenvalue of $T_m(u)$ whenever $m\equiv1\pmod4$.
\end{proposition}

\begin{proof}
Let $f_m:=D_{(1,\dots,1)}$. By \cite[Lemma~2.4]{Blattner} we have $D_s\equiv f_m\pmod2$ and $f_m(y+y^2)=(1+y)^{m+1}+y^{m+1}$ over $\FF_2$. Writing $m+1=2^aq$ with $q$ odd, the latter is the $2^a$-th power of $(1+y)^q+y^q$, and so has degree $2^a(q-1)$ in $y$. Thus the reduction of $f_m$ modulo $2$ has degree $(m+1-2^a)/2$, which gives $\deg D_s\ge(m+1-2^a)/2$ and hence the stated bound.

To see that equality can hold we induct on $a$, the cases $a=0$ and $m=1$ being trivial. If $a\ge1$ and $m=2n+1\ge3$, then $v_2(n+1)=a-1$, and by induction there is a sequence $\sigma$ of length $n-1$ with $\mu(\sigma)=2^{a-1}-1$. For $s:=(\sigma,\varepsilon,-\varepsilon,\rev\sigma)$ we have $D_s=D_\sigma^2$ by \cite[Lemma~2.3]{Blattner} (see also Lemma~\ref{lem:fold}), and therefore $\mu(s)=2\mu(\sigma)+1=2^a-1$.
\end{proof}

The bound shows in particular that $\deg D_s=1$ is possible only if $m=2$ or $m=2^K-1$. Our second result describes the latter case completely. For a sequence $\sigma$ we write $D_\sigma^-$ for the matching polynomial of the path of $\sigma$ with its last vertex removed; thus $D_\sigma^-=D_{k-1}$ if $D_\sigma=D_k$ in the recursion above.

\begin{proposition}\label{prop:deg1}
Let $m=2^K-1$ with $K\ge2$, and let $s\in\{-1,1\}^{m-1}$. Then $\deg D_s=1$ if and only if
\begin{equation}\label{eq:antifold}
s=(\sigma,\varepsilon,\varepsilon,-\rev\sigma)\qquad\text{for some }\sigma\in F_{K-1}\text{ and }\varepsilon\in\{-1,1\}.
\end{equation}
In this case $D_s=1+2\varepsilon x$. Thus there are exactly $2^{K-1}$ such sequences, and the corresponding matrices $T_m(u)$ have the eigenvalue zero with multiplicity $m-2$, together with the two eigenvalues $\pm\sqrt{-2\varepsilon}$.
\end{proposition}

\begin{proof}
Let $h:=(m-1)/2$, and let $c_e:=xf_{e-1}f_{m-1-e}$ for $1\le e\le m-1$, so that $c_e=c_{m-e}$. By \cite[Eq.~(5)]{Blattner} we have $D_s\equiv f_m-2\sum_et_ec_e\pmod4$, where $s_e=1-2t_e$. Comparing this with the same congruence for an element of $F_K$, for which $D=1$, we obtain
\begin{equation}\label{eq:rho}
\frac{D_s-1}{2}\equiv\sum_{e=1}^{h}\rho_e\,c_e\pmod2,
\end{equation}
where $\rho_e\in\{0,1\}$ records whether $s$ deviates from the pattern of a fold at the mirror pair $(e,m-e)$: we set $\rho_e:=1$ if $s_e\ne s_{m-e}$ and $e<h$, or if $s_h=s_{h+1}$ and $e=h$, and $\rho_e:=0$ otherwise. For $s=(1,\dots,1)$ this reads $(f_m-1)/2\equiv c_h\pmod2$, and from the formula for $f_m(y+y^2)$ one computes $c_h\equiv x+x^2+x^4+\dots+x^{2^{K-2}}\pmod2$. On the other hand, each matching with $k$ edges is counted $k$ times in $\sum_ec_e$, so that $2\sum_{e\le h}c_e=xf_m'(x)$. Combining these facts,
\[
\sum_{e=1}^{h}c_e=x\Bigl(\frac{f_m-1}{2}\Bigr)'\equiv x\,c_h'(x)\equiv x\pmod2.
\]
Since the $c_e$ with $e\le h$ are linearly independent modulo $2$ \cite[Lemma~3.1]{Blattner}, we see from \eqref{eq:rho} that $\deg D_s\le1$ forces $\rho=0$ or $\rho=(1,\dots,1)$.

If $\rho=0$, then $s=(\sigma,\varepsilon,-\varepsilon,\rev\sigma)$ and $D_s=D_\sigma^2$ has even degree, so that $D_s=1$. If $\rho=(1,\dots,1)$, then $s=(\sigma,\varepsilon,\varepsilon,-\rev\sigma)$ for some $\sigma$ and $\varepsilon$. As negating a sequence amounts to replacing $x$ by $-x$, expansion at the central vertex gives
\begin{equation}\label{eq:central}
D_s(x)=D_\sigma(x)D_\sigma(-x)+\varepsilon x\bigl(D_\sigma^-(x)D_\sigma(-x)+D_\sigma(x)D_\sigma^-(-x)\bigr).
\end{equation}
The first term on the right is the even part of $D_s$, and has degree $2\deg D_\sigma$. Thus $\deg D_s\le1$ forces $D_\sigma=1$, that is, $\sigma\in F_{K-1}$. In that case $D_\sigma^-=1\pm x\pm x^3\pm\dots\pm x^{2^{K-2}-1}$ by \cite[Corollary~4.3]{Blattner} applied to $\rev\sigma$ (see also Proposition~\ref{prop:resolvent}), so that the even part of $D_\sigma^-$ is $1$ and \eqref{eq:central} becomes $D_s=1+2\varepsilon x$. Conversely, this computation shows that every sequence of the form \eqref{eq:antifold} has $D_s=1+2\varepsilon x$.
\end{proof}

Table~\ref{tab:deg} lists the number of sequences $s$ with $\deg D_s=k$ for $k\le7$. For $m=15$ and $m=31$ it was obtained by enumerating all sequences, and for $m=63$ by enumerating the $2^{34}$ sequences which \eqref{eq:rho} permits when $\deg D_s\le3$.

\begin{table}[H]
\centering
\small
\begin{tabular}{rrrrrrrrr}
\toprule
$\deg D_s$ & 0 & 1 & 2 & 3 & 4 & 5 & 6 & 7\\
\midrule
$m=15$ & 8 & 8 & 24 & 0 & 126 & 520 & 3794 & 11904\\
$m=31$ & 16 & 16 & 48 & 0 & 128 & 16 & 48 & 32\\
$m=63$ & 32 & 32 & 96 & 0 & & & & \\
\bottomrule
\end{tabular}
\caption{The number of sequences $s\in\{-1,1\}^{m-1}$ with $\deg D_s=k$, for $k\le7$. The eigenvalue zero of $T_m(u)$ then has multiplicity $m-2k$.}\label{tab:deg}
\end{table}

The first two columns are explained by \cite[Theorem~2.2]{Blattner} and by Proposition~\ref{prop:deg1}. The most striking feature of the table is the absence of the degree $3$. For $m=15$ and $m=31$ it is the only degree that does not occur, the degrees from $8$ to $15$ being present for $m=31$ as well; for $m=15$ this is the observation of Corless mentioned above. We are led to the following conjecture.

\begin{conjecture}\label{conj:deg3}
Let $m=2^K-1$ with $K\ge4$. Then there are exactly $3\cdot2^{K-1}$ sequences $s\in\{-1,1\}^{m-1}$ with $\deg D_s=2$, and none with $\deg D_s=3$. In particular, the eigenvalue zero of $T_m(u)$ never has multiplicity $m-6$.
\end{conjecture}

\begin{remark}\label{rem:twists}
The argument of Proposition~\ref{prop:deg1} goes some way towards this conjecture. If $\deg D_s\le3$, then \eqref{eq:rho} confines $\rho$ to a three-dimensional subspace of $\FF_2^h$, which contains $(1,\dots,1)$. Two of its eight elements can be dealt with for all $K\ge4$. If $\rho=0$, then $D_s=D_\sigma^2$ with $\deg D_\sigma\le1$, so that $D_s=1$ or $D_s=(1\pm2x)^2$ by Proposition~\ref{prop:deg1}. If $\rho=(1,\dots,1)$, then the even part of \eqref{eq:central} again gives $\deg D_\sigma\le1$. The case $D_\sigma=1$ leads to $D_s=1\pm2x$ as before, while the case $\deg D_\sigma=1$ cannot occur. Indeed, in that case $D_\sigma=1+\beta x$ with $\beta=\pm2$ by Proposition~\ref{prop:deg1}, and $D_\sigma^-$ has odd degree $L:=2^{K-2}-1$ and leading coefficient $\pm1$, because the path of $\sigma$ with its last vertex removed has a unique perfect matching. Writing $E$ and $O$ for the even and odd parts of $D_\sigma^-$, the odd part of \eqref{eq:central} is $2\varepsilon x(E(x)-\beta xO(x))$, which has degree $L+2\ge5$. The six remaining elements have been settled only by computation, for $K\le6$. For $m=31$ the polynomials $D_s$ of degree $2$ are $(1\pm2x)^2$, $1\pm2x^2$ and $1\pm2x\pm2x^2$.
\end{remark}

\enlargethispage{2\baselineskip}
\section{Formal verification and computations}\label{sec:lean}

We have formalised the proofs of Theorem~\ref{thm:main}, Corollary~\ref{cor:sizes}, Theorem~\ref{thm:mod4} and Corollary~\ref{cor:power}(i) in Lean~4 \cite{Lean}, version 4.36.0, using the Mathlib library \cite{mathlib}. The formalisation follows Sections~\ref{sec:prelim} and~\ref{sec:proofs} and has about $3{,}000$ lines. The characteristic polynomial of a pair is defined by its three-term recurrence. The tridiagonal matrix with diagonal $d$, superdiagonal $a$ and subdiagonal $c$, with rows and columns indexed from $0$, is defined by
\begin{Verbatim}
def Tri (n : ℕ) (d a c : ℕ → R) : Matrix (Fin n) (Fin n) R :=
  fun i j =>
    if (i : ℕ) = j then d i
    else if (j : ℕ) = i + 1 then a i
    else if (i : ℕ) = j + 1 then c j
    else 0
\end{Verbatim}
\noindent and its characteristic polynomial is shown to be that of the pair $(d,ac)$. The pairs obtained from $o$ by folds are described by the inductive predicate
\begin{Verbatim}
inductive Folded : List ℤ → List ℤ → Prop
  | base : Folded [0] []
  | odd {d b : List ℤ} {ε : ℤ} (hp : Folded d b) (hε : ε = 1 ∨ ε = -1) :
      Folded (d ++ 0 :: d.reverse) (b ++ ε :: -ε :: b.reverse)
  | even {d b : List ℤ} {δ : ℤ} (hp : Folded (d ++ [0]) b) (hδ : δ = 1 ∨ δ = -1) :
      Folded (d ++ -δ :: δ :: d.reverse) (b ++ -1 :: b.reverse)
\end{Verbatim}
\noindent in which the even fold is applied to the pairs whose last diagonal entry is $0$. The formal statement of Theorem~\ref{thm:main} is
\begin{Verbatim}
theorem isNilpotent_iff_folded (d a c : List ℤ) (ha : a.length + 1 = d.length)
    (hc : c.length + 1 = d.length) (hd : ∀ x ∈ d, x = -1 ∨ x = 0 ∨ x = 1)
    (ha1 : ∀ x ∈ a, x = 1 ∨ x = -1) (hc1 : ∀ x ∈ c, x = 1 ∨ x = -1) :
    IsNilpotent (Tri d.length (seqOf d) (seqOf a) (seqOf c))
      ↔ Folded d (List.zipWith (· * ·) a c)
\end{Verbatim}
\noindent where \texttt{seqOf} extends a list by zeros to a sequence; the uniqueness assertion of the theorem is the injectivity of the two folds. With \mbox{\texttt{NilPairs n}} denoting the set $\cN_n$, Corollary~\ref{cor:sizes} reads
\begin{Verbatim}
theorem nilPairs_nonempty_iff (n : ℕ) :
    (NilPairs n).Nonempty ↔ n = 1 ∨ ∃ j : ℕ, 3 * 2 ^ j - 1 ≤ n ∧ n ≤ 2 ^ (j + 2) - 1
theorem nilPairs_ncard {n : ℕ} (hn : (NilPairs n).Nonempty) :
    (NilPairs n).ncard = 2 ^ Nat.log 2 n
\end{Verbatim}
\noindent For a pair given by two sequences $d$ and $b$ of integers, with characteristic polynomial \mbox{\texttt{K d b n}}, and with \mbox{\texttt{InC d b n}} expressing that the pair lies in $\cC_n$, Corollary~\ref{cor:power}(i) reads
\begin{Verbatim}
theorem nil_of_congr (n : ℕ) : ∀ {d b : ℕ → ℤ}, 1 ≤ n → InC d b n →
    (2 : ℤ[X]) ^ (Nat.log 2 n + 1) ∣ K d b n - X ^ n → K d b n = X ^ n
\end{Verbatim}
\noindent Theorem~\ref{thm:mod4} is stated in the same setting, with the property of being a fold expressed by equations between the entries, and with \eqref{eq:descent} in the form that for every $k\ge1$ the polynomial $\chi_P-\lambda^n$ is divisible by $2^{k+1}$ if and only if $\chi_p-\lambda^h$ is divisible by $2^k$. All these theorems depend only on the three standard axioms of Lean (\texttt{propext}, \texttt{Classical.choice} and \texttt{Quot.sound}). Corollary~\ref{cor:count}, parts (ii) and (iii) of Corollary~\ref{cor:power}, Proposition~\ref{prop:resolvent} and the results of Sections~\ref{sec:liu} and~\ref{sec:mult} have not been formalised.

The results of this paper also lend themselves to verification by computer. We have confirmed Theorem~\ref{thm:main} by an exhaustive search for all $n\le27$, Theorem~\ref{thm:mod4} and Corollary~\ref{cor:power} for all pairs with $n\le10$, the values in Corollary~\ref{cor:power}(ii) and (iii) for $n\le21$, and Corollary~\ref{cor:count} by enumeration of all matrices for $n\le6$. The formalisation and the programs, which also produce Table~\ref{tab:deg} and the values of $J(n)$ in Remark~\ref{rem:full}, are available at \url{https://github.com/marcelbtec/nilpotent-tridiagonal} and are archived at \cite{BlattnerCode}.

\section{Some questions}\label{sec:questions}
We close with some questions suggested by the above results.
\begin{enumerate}
\item Remark~\ref{rem:population} shows that Theorem~\ref{thm:main} is sensitive to the population. Is there a similar description of the nilpotent pairs $(d,b)$ with $d\in\ZZ^n$ and $b\in\{-1,1\}^{n-1}$, or of those in which the $b_k$ are allowed to be arbitrary complex numbers of modulus $1$?
\item Prove Conjecture~\ref{conj:deg3}. More generally, which degrees of $D_s$ occur for a given $m$?
\item Liu's matrices of sizes $6$, $14$ and $28$ are folds, while the one of size $5$ is not. For which dimensions does Liu's construction produce entries in $\{-1,0,1\}$?
\end{enumerate}

\subsection*{Acknowledgements}
The formulation of Theorem~\ref{thm:mod4} as a separate statement, parts (i) and (ii) of Corollary~\ref{cor:power}, the test of Remark~\ref{rem:algo}, the example and the lower bound for $J(7)$ in Remark~\ref{rem:full}, the second example in Remark~\ref{rem:population}, the reference to the instance of size $5$ in \cite{Berhanu}, and a gap in an earlier version of Remark~\ref{rem:twists} were pointed out in reports on earlier drafts of this paper.

\subsection*{Collaboration with an AI system}
This paper is the result of a collaboration between the author and the AI system Claude (Anthropic; Fable~5.2). The system wrote the Lean formalisation and the programs of Section~\ref{sec:lean}, carried out the literature search. The author posed the questions, worked out proof sketches, and directed the work. The author has checked all arguments and takes full responsibility for the paper.

\end{document}